\documentclass[a4paper,12pt]{amsart}
\usepackage{amsfonts}
\usepackage{amssymb}
\usepackage{ifthen}
\usepackage{amscd}
\usepackage{amsxtra}
\usepackage{graphicx}
\usepackage{color}
\nonstopmode \numberwithin{equation}{section}
\newtheorem{thm}{Theorem}
\newtheorem{lem}{Lemma}
\newtheorem{cor}{Corollary}
\newtheorem{prop}{Proposition}

\newtheorem{cl}{Claim}
\newtheorem{ca}{Case}
\newtheorem{sca}{Subcase}
\newtheorem{scl}{Subclaim}
\newtheorem{conj}{Conjecture}

\theoremstyle{definition}
\newtheorem{defn}{Definition}

\newtheorem{op}[equation]{Open Problem}
\newtheorem{ques}[equation]{Question}
\newtheorem{rem}{Remark}[section]
\newtheorem{exam}[equation]{Example}

\newcounter {own}
\def\theown {\thesection       .\arabic{own}}

\newenvironment{pf}[1][]{%
 \vskip 3mm
 \noindent
 \ifthenelse{\equal{#1}{}}%
  {{\slshape Proof. }}%
  {{\slshape #1.} }%
 }%
{\qed\bigskip}

\newcounter{alphabet}
\newcounter{tmp}
\newenvironment{Thm}[1][]{\refstepcounter{alphabet}%
\bigskip%
\noindent%
{\bf Theorem \Alph{alphabet}}%
\ifthenelse{\equal{#1}{}}{}{ (#1)}%
{\bf .} \itshape}{\vskip 8pt}

\makeatletter
\newcommand{\Ref}[1]{\@ifundefined{r@#1}{}{\setcounter{tmp}{\ref{#1}}\Alph{tmp}}}
\makeatother

\newenvironment{Lem}[1][]{\refstepcounter{alphabet}%
\bigskip%
\noindent%
{\bf Lemma \Alph{alphabet}}%
{\bf .} \itshape}{\vskip 8pt}

\def\be{\begin{equation}}
\def\ee{\end{equation}}

\newcommand{\bee}{\begin{enumerate}}
\newcommand{\eee}{\end{enumerate}}

\newcommand{\blem}{\begin{lem}}
\newcommand{\elem}{\end{lem}}
\newcommand{\bthm}{\begin{thm}}
\newcommand{\ethm}{\end{thm}}
\newcommand{\bcor}{\begin{cor}}
\newcommand{\ecor}{\end{cor}}
\newcommand{\beg}{\begin{exam}}
\newcommand{\eeg}{\end{exam}}
\newcommand{\begs}{\begin{examples}}
\newcommand{\eegs}{\end{examples}}
\newcommand{\bdefe}{\begin{defn}}
\newcommand{\edefe}{\end{defn}}
\newcommand{\bprob}{\begin{prob}}
\newcommand{\eprob}{\end{prob}}
\newcommand{\bques}{\begin{ques}}
\newcommand{\eques}{\end{ques}}
\newcommand{\bei}{\begin{itemize}}
\newcommand{\eei}{\end{itemize}}
\newcommand{\bcon}{\begin{conj}}
\newcommand{\econ}{\end{conj}}
\newcommand{\bop}{\begin{op}}
\newcommand{\eop}{\end{op}}

\newcommand{\bca}{\begin{ca}}
\newcommand{\eca}{\end{ca}}
\newcommand{\bsca}{\begin{sca}}
\newcommand{\esca}{\end{sca}}

\newcommand{\bcl}{\begin{cl}}
\newcommand{\ecl}{\end{cl}}

\newcommand{\bscl}{\begin{scl}}
\newcommand{\escl}{\end{scl}}

\newcommand{\bcons}{\begin{conjs}}
\newcommand{\econs}{\end{conjs}}
\newcommand{\bprop}{\begin{propo}}
\newcommand{\eprop}{\end{propo}}
\newcommand{\br}{\begin{rem}}
\newcommand{\er}{\end{rem}}
\newcommand{\brs}{\begin{rems}}
\newcommand{\ers}{\end{rems}}
\newcommand{\bo}{\begin{obser}}
\newcommand{\eo}{\end{obser}}
\newcommand{\bos}{\begin{obsers}}
\newcommand{\eos}{\end{obsers}}
\newcommand{\bpf}{\begin{pf}}
\newcommand{\epf}{\end{pf}}
\newcommand{\ba}{\begin{array}}
\newcommand{\ea}{\end{array}}
\newcommand{\beq}{\begin{eqnarray}}
\newcommand{\beqq}{\begin{eqnarray*}}
\newcommand{\eeq}{\end{eqnarray}}
\newcommand{\eeqq}{\end{eqnarray*}}

\newcommand{\ds}{\displaystyle}

\newcounter{minutes}
\divide\time by 60
\newcounter{hours}
\multiply\time by 60 \addtocounter{minutes}{-\time}

\begin{document}

\bibliographystyle{amsplain}
\title []
{Lipschitz type, radial growth and Dirichlet type spaces on
functions induced by certain elliptic
 operators}

\def\thefootnote{}
\footnotetext{ \texttt{\tiny File:~\jobname .tex,
          printed: \number\day-\number\month-\number\year,
          \thehours.\ifnum\theminutes<10{0}\fi\theminutes}
} \makeatletter\def\thefootnote{\@arabic\c@footnote}\makeatother

\author{Shaolin Chen}
\address{S. Chen, College of Mathematics and
Statistics, Hengyang Normal University, Hengyang, Hunan 421008,
People's Republic of China.}  \email{mathechen@126.com}





\author{Antti Rasila }
\address{A. Rasila, Department of Mathematics and Systems Analysis, Aalto University, P. O. Box 11100, FI-00076 Aalto,
 Finland.} \email{antti.rasila@iki.fi}


\begin{abstract} In this paper, we  investigate properties of
classes of functions related to certain elliptic
 operators. Firstly, we prove that a main result of Dyakonov (Acta Math. 178(1997),
 143--167) on analytic functions can be extended to this more general
 setting. Secondly, we study the radial growth on these functions and the
 obtained results are  generalizations of  the corresponding results
 of Makarov (Proc. London Math. Soc. 51(1985), 369--384) and Korenblum (Bull. Amer. Math. Soc. 12(1985), 99--102).
Finally, we discuss the Dirichlet type energy integrals on such
classes of functions and their applications.
\end{abstract}



\maketitle \pagestyle{myheadings} \markboth{ S. Chen and A.
Rasila}{Lipschitz type, radial growth and Dirichlet type spaces}

\section{Introduction and  main results }\label{csw-sec1}

Let $\mathbb{R}^{n}
$ denote the $n$-dimensional Euclidean space, where
$n\geq2$. For $a=(a_{1},\ldots,a_{n})$ and
$x=(x_{1},\ldots,x_{n})\in\mathbb{R}^{n}$, we define the
inner product $\langle \cdot ,\cdot \rangle$ by $\langle
x,a\rangle=x_{1}a_{1}+\cdots+x_{n}a_{n} $ so that the Euclidean
length of $x$ is defined by $|x|=\langle
x,x\rangle^{1/2}=(|x_{1}|^{2}+\cdots+|x_{n}|^{2})^{1/2}. $ Denote a
ball in $\mathbb{R}^{n}$ with center $x'$ and radius $r$ by
$\mathbb{B}^{n}(x',r)=\{x\in\mathbb{R}^{n}:\, |x-x'|<r\}. $ In
particular, let $\mathbb{B}^{n}=\mathbb{B}^{n}(0,1)$ and
$\mathbb{B}_{r}^{n}=\mathbb{B}^{n}(0,r)$. Set
$\mathbb{D}=\mathbb{B}^2$, the open unit disk in the complex plane
$\mathbb{C}$. Let $\Omega$ be a domain of $\mathbb{R}^{n}$, with
non-empty boundary. We use $d_{\Omega}(x)$ to denote the Euclidean
distance from $x\in \Omega$ to the boundary $\partial \Omega$ of $\Omega$.
If $\Omega=\mathbb{B}^{n}$, we write $d(x)$ instead of  $d_{\mathbb{B}^{n}}(x).$
We denote by $ \mathcal{C}^{m}(\Omega)$ the set of all $m$-times
continuously differentiable  functions from a domain
$\Omega\subset\mathbb{R}^{n}$
 into $\mathbb{R}$, where $m\in\{1,2,\ldots\}$. Furthermore,  we use $C$ to denote the various positive
constants, whose value may change from one occurrence to another.

Fix $\tau\geq1$, and let $u\in
\mathcal{C}^{2}(\mathbb{B}^{n})$ be a solution to the equation
\be\label{eq1c}\Delta u=\lambda(x)|u|^{\tau-1}u,\ee where $\Delta$ is
the Laplacian operator and
 $\lambda$ is a  continuous
function from $\mathbb{B}^{n}$ into $\mathbb{R}$. If $\lambda$ is a
constant function in (\ref{eq1c}), then this type equation has
attracted the attention of many authors, where the case $\tau=1$ and $\lambda<0$,
i.e., the {\it Helmholtz} equation, is particularly important. We refer to \cite{Ba,
Cl,E,Fa, Pas} and the references therein.
If $\lambda>0$ is a constant and $\tau=1$, then (\ref{eq1c}) is the
{\it Yukawa equation}, which arose out of an attempt of the Japanese Nobel physicist Hideki Yukawa to describe the
nuclear potential of a point charge as $e^{-\sqrt{\lambda} r}/r$
(cf. \cite{A,BS,CPR,CRV,CRW,D-,D-1, SW,St,Ya}). It is well known
that if $\lambda$ is a constant function and $\tau=1$, then each
solution $u$ to (\ref{eq1c}) belongs to
$\mathcal{C}^{\infty}(\mathbb{B}^{n}).$ Moreover, if $\lambda=0$ in
(\ref{eq1c}), then $u$ is harmonic in $\mathbb{B}^{n}$.

In fact, the equation (\ref{eq1c}) can be regarded as the induced
equation by the elliptic partial differential operators $\mbox{div}
p^{2}\nabla +q$, where $\nabla$ denotes the gradient and  $p$,  $q$
are real-valued functions satisfying $p\in
\mathcal{C}^{2}(\mathbb{B}^{n})$ and $p\neq0$ in $\mathbb{B}^{n}$.
Precisely, the elliptic operators \be\label{eq-1}E_{p,q}=\mbox{div}\,
p^{2}\nabla +q\ee can be decomposed into the following form (cf.
\cite{N,V}) \be\label{eq-1.0}E_{p,q}=p\big(\Delta-\varphi\big)p,\ee
where
  $\varphi=(\Delta p)/p-q/p^{2}.$ By
(\ref{eq-1.0}), we see that the equation
\be\label{eq-2}E_{p,q}(u)=\big(\mbox{div}\, p^{2}\nabla
+q\big)u=0~\mbox{in}~\mathbb{B}^{n}\ee  is equivalent to {\it the
stationary Schr\"odinger type equation} (cf. \cite{A, V})
\be\label{eq-1.1}\Delta h=\varphi h,\ee where $h=pu.$  If we can
choose  some $p$ and $q$ such that $\varphi=\lambda |h|^{\tau-1}$,
then (\ref{eq-1.1}) is the same type equation as (\ref{eq1c}), where
$\tau\geq1$. In particular, if $n=2$, the equation (\ref{eq-2})
 is closely related to {\it the main Vekua equation} (cf.
\cite{B1, B2, V,Ve})

\be\label{eq-3}
\partial_{\overline{z}}w=\frac{\partial_{\overline{z}}f}{f}\overline{w},\ee where $z=x+iy$, $\partial_{z}=\frac{1}{2}
(\partial/\partial x-i\partial/\partial y)$ and
$\partial_{\overline{z}}=\frac{1}{2} (\partial/\partial
x+i\partial/\partial y)$. In fact, if $f=pu_{0}$, then, for any
solution $u$ to the equation (\ref{eq-2}), there is is a
corresponding solution $w$ to the equation (\ref{eq-3}) such that
$u=\mbox{Re}w/p$ is a solution to the equation (\ref{eq-2}), where
$u_{0}$ is a positive solution to the equation (\ref{eq-2}).


\begin{prop}\label{prop} Suppose
 $u\in \mathcal{C}^{2}(\mathbb{B}^{n})$ is  a solution to the equation {\rm(\ref{eq1c})}, where
 $\lambda$ is a nonnegative continuous
function from $\mathbb{B}^{n}$ into $\mathbb{R}$ with
$\sup_{x\in\mathbb{B}^{n}}\lambda(x)<+\infty$.
 For all $x\in\mathbb{B}^{n},$ there is a
positive constant $C$  such that
$$|\nabla
u(x)|^{\nu}\leq\frac{C}{R^{\nu+n}}\left(\int_{\mathbb{B}^{n}(x,R)}|u(y)|^{\nu}dy\\
+\int_{\mathbb{B}^{n}(x,R)}|u(y)|^{\tau\nu}dy\right),
$$ where  $\nu\in[1,+\infty)$  and $R$ is a positive constant such that
$\overline{\mathbb{B}^{n}(x,R)}\subset\mathbb{B}^{n}$.
\end{prop}

A continuous increasing function $\omega:\, [0,+\infty)\rightarrow
[0,+\infty)$ with $\omega(0)=0$ is called a {\it majorant} if
$\omega(t)/t$ is non-increasing for $t>0$ (cf.
\cite{CP,CPR,CPR-2015,Dy1,Dy2,P}). Given a subset $\Omega$ of
$\mathbb{R}^{n}$, a function $u:\, \Omega\rightarrow \mathbb{R}$ is
said to belong to the {\it Lipschitz space $L_{\omega}(\Omega)$} if
there is a
positive constant $C$ such that 
$$|u(x_{1})-u(x_{2})|\leq C\omega(|x_{1}-x_{2}|) ~\mbox{ for all
$x_{1},x_{2}\in\Omega.$}$$

For $\nu\in(0,+\infty]$, the {\it generalized Hardy space
$\mathcal{H}^{\nu}_{g}(\mathbb{B}^{n})$} consists of all those
functions $u:\mathbb{B}^{n}\rightarrow\mathbb{R}$ such that $u$ is
measurable, $M_{\nu}(r,f)$ exists for all $r\in(0,1)$ and  $
\|u\|_{\nu}<+\infty$, where
$$\|u\|_{\nu}=
\begin{cases}
\displaystyle\sup_{0<r<1}M_{\nu}(u,r)
& \mbox{if } \nu\in(0,+\infty),\\
\displaystyle\sup_{x\in\mathbb{B}^{n}}|u(x)| &\mbox{if }
\nu=+\infty,
\end{cases}
~
M_{\nu}(u,r)=\left(\int_{\partial\mathbb{B}^{n}}|u(r\zeta)|^{\nu}\,d\sigma(\zeta)\right)^{\frac{1}{\nu}},
$$ and $d\sigma$ denotes the normalized surface measure on $\partial\mathbb{B}^{n}$.

The classical  {\it harmonic Hardy space
$\mathcal{H}^{p}(\mathbb{D})$} consisting of harmonic functions in
$\mathbb{D}$ is a subspace of $\mathcal{H}^{p}_{g}(\mathbb{D})$.

\begin{defn}\label{def-1}
For $\nu\in(0,+\infty]$, $\alpha>0$, $\beta\in\mathbb{R}$ and a majorant $\omega$, 
we use
$\mathcal{L}_{\nu,\omega}\mathcal{B}^{\beta}_{\alpha}(\mathbb{B}^{n})$
to denote the {\it generalized Bloch type space} of all  functions
$u\in \mathcal{C}^{1}(\mathbb{B}^{n})$  with
$\|u\|_{\mathcal{L}_{\nu,\omega}\mathcal{B}^{\beta}_{\alpha}(\mathbb{B}^{n})}\\
<+\infty$, where
$$\|u\|_{\mathcal{L}_{\nu,\omega}\mathcal{B}^{\beta}_{\alpha}(\mathbb{B}^{n})}=
\begin{cases}
\displaystyle|u(0)|+\sup_{x\in\mathbb{B}^{n}} \left\{
M_{\nu}(|\nabla u|,|x|)\omega\big(\phi(x)\big)\right\}
& \mbox{if } \nu\in(0,+\infty),\\
\displaystyle|u(0)|+\sup_{x\in\mathbb{B}^{n}} \left\{ |\nabla
u(x)|\omega\big(\phi(x)\big)\right\} &\mbox{if } \nu=+\infty,
\end{cases}
$$  and $\phi(x)=d^{\alpha}(x)(1-\log\ d(x))^{\beta}$.
\end{defn}

It is easy to see  that
$\mathcal{L}_{\nu,\omega}\mathcal{B}^{\beta}_{\alpha}(\mathbb{B}^{n})$
is a Banach space for $\nu\geq1$. Moreover, we have the following:

\begin{enumerate}
\item[{\rm (1)}] If  $\beta=0$, then
$\mathcal{L}_{+\infty,\omega}\mathcal{B}^{0}_{\alpha}(\mathbb{D})$
is called the {\it $\omega$-$\alpha$-Bloch space} (cf. \cite{CPR}).

\item[{\rm (2)}] If we take $\alpha=1$, then
$\mathcal{L}_{+\infty,\omega}\mathcal{B}^{\beta}_{1}(\mathbb{D})$ is
called  the {\it  logarithmic $\omega$-Bloch space}.

\item[{\rm (3)}] If we take $\omega(t)=t$ and $\beta=0$,
then
$\mathcal{L}_{+\infty,\omega}\mathcal{B}^{0}_{\alpha}(\mathbb{D})$
is called  the  {\it generalized $\alpha$-Bloch space} (cf.
\cite{CPR, RU,Z,Z1}).

\item[{\rm (4)}] If we take $\omega(t)=t$ and $\alpha=1$, then
$\mathcal{L}_{+\infty,\omega}\mathcal{B}^{\beta}_{1}(\mathbb{D})$ is
called  the {\it generalized logarithmic Bloch space} (cf.
\cite{CPR, CPW1,Dyak,GPPJ,Pav1,Pe,Z1}).
\end{enumerate}

In \cite{Kr}, the author studied the Lipschitz spaces on smooth
functions. Dyakonov \cite{Dy1} discussed the relationship between
the Lipschitz space and the bounded mean oscillation on analytic
functions in $\mathbb{D}$, and obtained the following result.

\begin{Thm}{\rm (\cite[Theorem 1]{Dy1})}\label{Thm-Day}
Suppose that $f$ is a analytic function in $\mathbb{D}$ which is
continuous up to the boundary of $\mathbb{D}$. If $\omega$ and
$\omega^{2}$ are regular majorants,  then
$$f\in L_{\omega}(\mathbb{D})\Longleftrightarrow \left(\mathcal{P}_{|f|^{2}}(z)-|f(z)|^{2}\right)^{1/2}\leq C\omega(d(z)),$$
where
\[
\mathcal{P}_{|f|^{2}}(z)=\frac{1}{2\pi}\int_{0}^{2\pi}\frac{1-|z|^{2}}{|z-e^{i\theta}|^{2}}|f(e^{i\theta})|^{2}d\theta,
\]
and $C$ is a positive constant.
\end{Thm}

In \cite{CPR,CRW}, the authors extended Theorem \Ref{Thm-Day} to
complex-valued harmonic functions (see \cite[Theorem 4]{CPR} and
\cite[Theorem 3]{CRW}). For the solutions to (\ref{eq1c}), we  get
the following  result, which is a generalization of Theorem
\Ref{Thm-Day}, \cite[Theorem 4]{CPR} and \cite[Theorem 3]{CRW}.

\begin{thm}\label{01-thm}
Let $\alpha\in[1,2)$ and  $\omega$ be a majorant.  Suppose that $u$
is a solution to {\rm (\ref{eq1c})} with $\tau=1$, where
 $\lambda$ is a nonnegative constant. Then
$u\in\mathcal{L}_{+\infty,\omega}\mathcal{B}^{0}_{\alpha}(\mathbb{B}^{n})$
if and only if there is a positive constant $C$ such that, for all
$r\in(0,d(x)]$,

\be\label{eq-csw13}\frac{1}{|\mathbb{B}^{n}(x,r)|}\int_{\mathbb{B}^{n}(x,r)}|u(y)-u(x)|dy\leq\frac{Cr}{\omega(r^{\alpha})},\ee
where $|\mathbb{B}^{n}(x,r)|$ denotes the volume of
$\mathbb{B}^{n}(x,r)$.
\end{thm}

Let $\Omega$ be a proper subdomain of $\mathbb{R}^{n}$. For
$x,y\in\Omega$, let
$$r_{\Omega}(x,y)=\frac{|x-y|}{\min\{d_{\Omega}(x),
d_{\Omega}(y)\}}.$$ The distance ratio metric (see e.g. \cite{Vu}) is defined for
$x,y\in\Omega$ by setting
$$j_{\Omega}(x,y)=\log(1+r_{\Omega}(x,y)).$$
We say that $f:~\Omega\rightarrow f(\Omega)\subset\mathbb{R}^{n}$
is  {\it weakly uniformly bounded} in $\Omega$
(with respect to $r_{\Omega}$)  if there is a constant $C>0$ such
that $r_{\Omega}(x,y)\leq1/2$ implies $r_{f(\Omega)}(f(x), f(y))\leq
C.$ For $x,y\in\Omega$, let
$$k_{\Omega}(x, y)=\inf_{\gamma}\int_{\gamma}\frac{ds}{d_{\Omega}(x)},$$
where infimum is taken over all rectifiable arcs
$\gamma\subset\Omega$ and $ds$ stands for the arc length measure on
$\gamma$ (cf. \cite{MM,Vu}).

In \cite{MM}, Mateljevi\'c and Vuorinen proved the following result.

\begin{Thm}{\rm (\cite[Theorem 2.8]{MM})} \label{Thm-MM}
Suppose that $\Omega$ is a proper subdomain of $\mathbb{R}^{n}$ and
$h:~\Omega\rightarrow\mathbb{R}^{n}$ is a harmonic mapping. Then the
following conditions are equivalent.

\begin{enumerate}
\item[{\rm (a)}] $h$ is weakly uniformly bounded;

\item[{\rm (b)}] There exists a constant $C$ such that, for all $x, y\in G$,
$$k_{u(\Omega)}(u(x),u(y))\leq Ck_{\Omega}(x,y).
$$
\end{enumerate}

\end{Thm}

We extended Theorem \Ref{Thm-MM} to the solutions of (\ref{eq1c})
with $\tau=1$, which is as follows.

\begin{thm}\label{02-thm}
Let $u=(u_{1},\cdots,u_{n})$ be a vector-valued  function from
$\mathbb{B}^{n}$ into the domain
$u(\mathbb{B}^{n})\subset\mathbb{R}^{n}$ satisfying $\Delta
u_{k}=\lambda_{k}u_{k},$ where $k\in\{1,\ldots, n\}$ and
 $\lambda_{k}$ is a nonnegative constant. Then the
following conditions are equivalent.

\begin{enumerate}
\item[{\rm (1)}] $u$ is weakly uniformly  bounded;

\item[{\rm (2)}] There exists a constant $C$ such that, for all $x, y\in\mathbb{B}^{n}$,
$$k_{u(\mathbb{B}^{n})}(u(x),u(y))\leq Ck_{\mathbb{B}^{n}}(x,y).
$$
\end{enumerate}
\end{thm}

We remark that we can replace $\mathbb{B}^{n}$ by some proper
domains $\Omega\subset\mathbb{R}^{n}$ in Theorem \ref{02-thm}.

 Makarov \cite{Ma} proved that if $f$ is analytic in $\mathbb{D}$ with $\mbox{Re}f\in\mathcal{L}_{+\infty,\omega}\mathcal{B}^{0}_{1}(\mathbb{D})$,
then there is a positive constant $C$ such that
\be\label{eq-csw1}\limsup_{r\rightarrow1-}\frac{|f(r\zeta)|}{\sqrt{\log\frac{1}{1-r}\log\log\log\frac{1}{1-r}}}\leq
C\|\mbox{Re}f\|_{\mathcal{L}_{+\infty,\omega}\mathcal{B}^{0}_{1}(\mathbb{D})}\ee
for almost every $\zeta\in\partial\mathbb{D}$, where $r\in[0,1)$ and
$\omega(t)=t$. In particular, Korenblum \cite{Ko} showed that if $u$
is a real harmonic function in $\mathbb{D}$ with
$u\in\mathcal{L}_{+\infty,\omega}\mathcal{B}^{0}_{1}(\mathbb{D}),$
then there is a positive constant $C$ such that
\be\label{eq-csw2}\limsup_{r\rightarrow1^{-}}\frac{|u(r\zeta)|}{\sqrt{\log\frac{1}{1-r}}\log\log\frac{1}{1-r}}\leq
C\|u\|_{\mathcal{L}_{+\infty,\omega}\mathcal{B}^{0}_{1}(\mathbb{D})}\ee
for almost every $\zeta\in\partial\mathbb{D}$, where $\omega(t)=t$.
For related investigations on the radial growth of Bloch type
functions, we refer to \cite{CM,FM,GPP,GP,GK,Po,ST}.




Analogously to (\ref{eq-csw1}) and (\ref{eq-csw2}), for
$\nu\in(0,+\infty)$ and for functions in
$u\in\mathcal{C}^{2}(\mathbb{B}^{n})$, satisfying a Bloch-type
condition, we prove the following result.


\begin{thm}\label{eq-y} Let $\omega$ be a majorant, $\nu\in[2,+\infty),$ $\alpha>0$ and
$\beta\leq\alpha$.   Suppose $u\in\mathcal{C}^{2}(\mathbb{B}^{n})$
satisfying $u\Delta u\geq0$ and $\big(|\nabla u|^{2}+u\Delta
u\big)\in\mathcal{L}_{\nu,\omega}\mathcal{B}^{\beta}_{\alpha}(\mathbb{B}^{n}).$
Then, for $n\geq3$ and $r\in[0,1)$,
$$M_{\nu}(u,r)\leq\left[|u(0)|^{2}+\frac{\nu(\nu-1)\||\nabla u|^{2}+u\Delta
u\|_{\mathcal{L}_{\nu,\omega}\mathcal{B}^{\beta}_{\alpha}(\mathbb{B}^{n})}r^{2}}{\omega(1)(n-2)}\int_{0}^{1}\frac{t(1-t^{n-2})}{\phi(tr)}dt\right]^{\frac{1}{2}},$$
where $\phi$ is defined as in definition {\rm\ref{def-1}}.

Moreover, for $n=2$ and $r\in[0,1)$,
$$M_{\nu}(u,r)\leq\left[|u(0)|^{2}+\frac{\nu(\nu-1)\||\nabla u|^{2}+u\Delta
u\|_{\mathcal{L}_{\nu,\omega}\mathcal{B}^{\beta}_{\alpha}(\mathbb{B}^{n})}r^{2}}{\omega(1)}\int_{0}^{1}\frac{t\log\frac{1}{t}}{\phi(tr)}dt\right]^{\frac{1}{2}}.$$
\end{thm}

\begin{defn}\label{def-2}
For $m\in\{2,3,\ldots\}$, we denote by
$\mathcal{HZ}_{m}(\mathbb{B}^{n})$ the class of all functions
$u\in\mathcal{C}^{m}(\mathbb{B}^{n})$ satisfying {\it Heinz's type
nonlinear differential inequality} (cf. \cite{CPR,HZ})
\be\label{eq-5} |\Delta u(x)|\leq a_{1}(x)|\nabla
u(x)|^{b_{1}}+a_{2}(x)|u(x)|^{b_{2}}+a_{3}(x)~ \mbox{ {\rm
for}}~x\in\mathbb{B}^{n}, \ee where $a_{k}~(k\in\{1,2,3\})$ are
real-valued nonnegative continuous functions in $\mathbb{B}^{n}$ and
$b_{j}~(j\in\{1,2\})$ are nonnegative constants.
\end{defn}

\begin{thm}\label{thm-1}
Let $\omega$ be a majorant, $\nu\in[2,+\infty),$ $\alpha>0$,
$\beta\leq\alpha$ and $u\in \mathcal{HZ}_{2}(\mathbb{B}^{n})\cap
\mathcal{L}_{\nu,\omega}\mathcal{B}^{\beta}_{\alpha}(\mathbb{B}^{n})$
satisfying $\sup_{x\in\mathbb{B}^{n}}a_{1}(x)<+\infty$,
$\sup_{x\in\mathbb{B}^{n}}a_{2}(x)<\frac{2n}{v},$
$\sup_{x\in\mathbb{B}^{n}}a_{3}(x)<+\infty$, $b_{1}\in[0,1]$ and
$b_{2}\in[0,1]$. If $n\geq3$ and $u\Delta u\geq0$, then, for
$r\in[0,1)$,

\begin{eqnarray*}
M_{\nu}( u,
r)&\leq&\Bigg[|u(0)|^{2}+\frac{\nu(\nu-1)}{(n-2)\omega^{2}(1)}\|u\|^{2}_{\mathcal{L}_{\nu,\omega}
\mathcal{B}^{\beta}_{\alpha}(\mathbb{B}^{n})}r^{2}\int_{0}^{1}\frac{t(1-t^{n-2})}{\phi^{2}(rt)}dt\\
&&+\frac{\nu
\sup_{x\in\mathbb{B}^{n}}a_{1}(x)}{(n-2)\omega^{b_{1}}(1)}\|u\|^{b_{1}}_{\mathcal{L}_{\nu,\omega}
\mathcal{B}^{\beta}_{\alpha}(\mathbb{B}^{n})}r^{2}M_{\nu}(u,r)\int_{0}^{1}\frac{t(1-t^{n-2})}{\phi^{b_{1}}(rt)}dt\\
&&+\frac{\nu
\sup_{x\in\mathbb{B}^{n}}a_{2}(x)}{2n}r^{2}M_{\nu}^{1+b_{2}}(
u, r)\\
&&+\frac{\nu \sup_{x\in\mathbb{B}^{n}}a_{3}(x)}{2n}r^{2}M_{\nu}( u,
r)\Bigg]^{\frac{1}{2}},
\end{eqnarray*}
where $\phi$ is defined as in Definition {\rm\ref{def-1}}.

 In
particular, if $n=2$ and $u\Delta u\geq0$, then, for $r\in[0,1)$,
\begin{eqnarray*}
M_{\nu}( u,
r)&\leq&\Bigg[|u(0)|^{2}+\frac{\nu(\nu-1)}{\omega^{2}(1)}\|u\|^{2}_{\mathcal{L}_{\nu,\omega}
\mathcal{B}^{\beta}_{\alpha}(\mathbb{D})}r^{2}\int_{0}^{1}\frac{t\log\frac{1}{t}}{\phi^{2}(rt)}dt\\
&&+\frac{\nu
\sup_{x\in\mathbb{D}}a_{1}(x)}{\omega^{b_{1}}(1)}\|u\|^{b_{1}}_{\mathcal{L}_{\nu,\omega}
\mathcal{B}^{\beta}_{\alpha}(\mathbb{D})}r^{2}M_{\nu}(u,r)\int_{0}^{1}\frac{t\log\frac{1}{t}}{\phi^{b_{1}}(rt)}dt\\
&&+\frac{\nu
\sup_{x\in\mathbb{D}}a_{2}(x)}{4}r^{2}M_{\nu}^{1+b_{2}}(
u, r)\\
&&+\frac{\nu \sup_{x\in\mathbb{D}}a_{3}(x)}{4}r^{2}M_{\nu}( u,
r)\Bigg]^{\frac{1}{2}}.
\end{eqnarray*}
\end{thm}

We remark that Theorem \ref{thm-1} is a generalization of
\cite[Theorem 1]{CPR}. As an application of Theorem \ref{thm-1}, we obtain the following
result.

\begin{cor}
Let $\omega$ be a majorant, $\nu\in[2,+\infty),$  $\alpha>0$ and
$\beta\leq\alpha$. Suppose that $u$ is a solution to {\rm
(\ref{eq1c})} with $\tau=1$, where
 $\lambda$ is a nonnegative continuous
function from $\mathbb{B}^{n}$ into $\mathbb{R}$ with
$\sup_{x\in\mathbb{B}^{n}}\lambda(x)<\frac{\nu}{2n}.$ If $n\geq3$
and $u\in
\mathcal{L}_{\nu,\omega}\mathcal{B}^{\beta}_{\alpha}(\mathbb{B}^{n})$,
then, for $r\in[0,1)$, $$M_{\nu}( u, r)\leq\frac{1}{C^{\ast}}
\left(|u(0)|^{2}+\frac{\nu(\nu-1)}{(n-2)\omega^{2}(1)}\|u\|^{2}_{\mathcal{L}_{\nu,\omega}
\mathcal{B}^{\beta}_{\alpha}(\mathbb{B}^{n})}r^{2}\int_{0}^{1}\frac{t(1-t^{n-2})}{\phi^{2}(rt)}dt\right)^{\frac{1}{2}},$$
where
$C^{\ast}=\left(1-\frac{r^{2}\nu}{2n}\sup_{x\in\mathbb{B}^{n}}\lambda(x)\right)^{\frac{1}{2}}$
and  $\phi$ is defined as in definition {\rm\ref{def-1}}.

 In
particular, if $n=2$ and  $u\in
\mathcal{L}_{\nu,\omega}\mathcal{B}^{\beta}_{\alpha}(\mathbb{B}^{n})$,
then, for $r\in[0,1)$, $$M_{\nu}( u, r)\leq\frac{1}{C^{\ast}}
\left(|u(0)|^{2}+\frac{\nu(\nu-1)}{\omega^{2}(1)}\|u\|^{2}_{\mathcal{L}_{\nu,\omega}
\mathcal{B}^{\beta}_{\alpha}(\mathbb{D})}r^{2}\int_{0}^{1}\frac{t\log\frac{1}{t}}{\phi^{2}(rt)}dt\right)^{\frac{1}{2}}.$$

\end{cor}

For $\alpha,\gamma,\mu\in\mathbb{R}$,
$$\mathcal{D}_{\nabla u}(\alpha,\gamma,\mu)=\int_{\mathbb{B}^{n}}(1-|x|^{2})^{\alpha}|\nabla u(x)|^{\gamma}
\left(\sum_{1\leq j,k\leq
n}u^{2}_{x_{j}x_{k}}(x)\right)^{\mu}dx<+\infty$$ is called a {\it
 Dirichlet type energy integral} of $u$ defined in
$\mathbb{B}^{n}$ (\cite{CPR, CRV, CRW,E, HHL, SH, ST,W, Ya}).


In \cite{CRW}, the authors investigated  certain properties on the above
Dirichlet type energy integral. In the following, we extend
\cite[Theorem 4]{CRW}  to a higher order form and give an
application.

\begin{thm}\label{thm-3} Let $u\in \mathcal{C}^{2}(\mathbb{B}^{n})$ be a solution to the
equation {\rm (\ref{eq1c})} with $\tau=1$. For  $\alpha>0$,
$\mu\in[1,n/2]$ and $\nu\in[2,+\infty)$,  if $\mathcal{D}_{\nabla
u}(\alpha,0,\mu)<+\infty$, then
$$\int_{\mathbb{B}^{n}}\big(d(x)\big)^{\beta\nu}\Delta \left(|\nabla
u(x)|^{\nu}\right)dx<+\infty,$$ where
$\beta=\frac{n+\alpha}{2\mu}-1.$
\end{thm}

We recall that a real function $f$ is said to have a {\it harmonic
majorant} if there is a positive harmonic function $F$ in
$\mathbb{B}^{n}$ such that, for all $x\in\mathbb{B}^{n}$,
$|f(x)|\leq F(x)$ (cf. \cite{CP,Du1,Nu, St, Ya}). Concerning
harmonic majorants, it is well known that a subharmonic function $u$
defined in $\mathbb{D}$ has a harmonic majorant if and only if
$\sup_{0<r<1}M_{1}(u,r)<+\infty$  (see \cite[Theorem 3.37]{Ho}). For
the solutions to (\ref{eq1c}), we have

\begin{thm}\label{thm-4}
Let $u\in \mathcal{C}^{2}(\mathbb{B}^{n})$ be a solution to the
equation {\rm (\ref{eq1c})} with $\tau=1$. Suppose that  $\alpha>0$,
$\mu\in[1,n/2]$ and $\nu\in[2,+\infty)$ satisfying
$\frac{n+\alpha}{2\mu}-1=\frac{1}{\nu}$. If $\mathcal{D}_{\nabla
u}(\alpha,0,\mu)<+\infty$, then $|\nabla
u|\in\mathcal{H}^{\nu}_{g}(\mathbb{B}^{n})$ and $|\nabla u|^{\nu}$
has a majorant.
\end{thm}

The proofs of Proposition \ref{prop}, and Theorems \ref{01-thm},
\ref{02-thm}, \ref{eq-y}, \ref{thm-1}, \ref{thm-3} and \ref{thm-4}
will be presented in Section \ref{csw-sec2}.

\section{Proofs of the  main results }\label{csw-sec2}

\begin{lem}\label{lem-1cw}
Let $u\in \mathcal{C}^{2}(\mathbb{B}^{n})$ with $u\Delta u\geq0$ in
$\mathbb{B}^{n}$. Then, for $\nu\geq1$, $|u|^{\nu}$ is subharmonic
in $\mathbb{B}^{n}$.
\end{lem}

\bpf Let $\mathcal{Z}_{u}=\{x\in\mathbb{B}^{n}:~u(x)=0\}$. Then
$\mathcal{Z}_{u}$ is a close set, which gives that
$\mathbb{B}^{n}\backslash\mathcal{Z}_{u}$ is an open set.  By
calculations, for $x\in\mathbb{B}^{n}\backslash\mathcal{Z}_{u}$, we
get \be\label{eq-2c}
\Delta(|u(x)|^{\nu})=\nu(\nu-1)|u(x)|^{\nu-2}|\nabla
u(x)|^{2}+\nu|u(x)|^{\nu-2}u(x)\Delta u(x) \geq0, \ee which implies
that $|u|^{\nu}$ is subharmonic in
 $\mathbb{B}^{n}$.
 \epf

 \begin{cor}\label{cor-1}
For some $\tau\geq1$, let $u\in \mathcal{C}^{2}(\mathbb{B}^{n})$ be
a solution to the equation {\rm (\ref{eq1c})},  where
 $\lambda$ is a nonnegative continuous
function in $\mathbb{B}^{n}$.
Then, for $\nu\geq1$, $|u|^{\nu}$ is subharmonic in
$\mathbb{B}^{n}$.
\end{cor}

In \cite{P-1994}, Pavlovi\'c proved the following result.

\begin{Lem}\label{lem-eqc-3}
Suppose that $\Omega$ is a bounded domain of $\mathbb{R}^{n}$ and
$u$ is a subharmonic function in $\Omega$.  For any $x\in\Omega$,
let $r$ be a positive constant such that
$\overline{\mathbb{B}^{n}(x,r)}\in\Omega$. Then, for $\nu>0$, there
are positive constant $C$  such that
$$|u(x)|^{\nu}\leq\frac{C}{r^{n}}\int_{\mathbb{B}^{n}(x,r)}|u(y)|^{\nu}dy.$$
\end{Lem}

The following result is well-known.

\begin{lem}\label{Lemx}
Suppose that $a,b\in[0,\infty)$ and $\iota\in(0,\infty)$. Then
$$(a+b)^{\iota}\leq2^{\max\{\iota-1,0\}}(a^{\iota}+b^{\iota}).
$$
\end{lem}

\subsection*{Proof of Proposition \ref{prop}}
 Let  $u\in \mathcal{C}^{2}(\mathbb{B}^{n})$ be a solution to
the equation (\ref{eq1c}).  Without loss of generality, we assume
that $x=0$ and  $n\geq3$. For $r\in(0,1)$ and all
$w\in\mathbb{B}^{n}_{r}$,

\be\label{eq-016csw}u(w)=r^{n-2}\left[\int_{\partial\mathbb{B}^{n}}P_{r}(w,\zeta)u(r\zeta)d\sigma(\zeta)-\int_{\mathbb{B}^{n}}G_{r}(w,y)\lambda(ry)|u(ry)|^{\tau-1}u(ry)dy\right],\ee
where $V(\mathbb{B}^{n})$ is the volume of the unit ball,
$$P_{r}(w,\zeta)=\frac{r^{2}-|w|^{2}}{|w-r\zeta|^{n}}$$ is the
Poisson kernel and
$$G_{r}(w,y)=\frac{1}{n(n-2)V(\mathbb{B}^{n})}\left[\frac{1}{|w-ry|^{n-2}}-\frac{1}{\big(r^{2}+|w|^{2}|y|^{2}-2r<w,y>\big)^{\frac{n-2}{2}}}\right]$$
is the Green function (see \cite{GT, Ho}). By calculations, we have

\be\label{eq-1c-2}|\nabla
P_{r}(0,\zeta)|=O\left(\frac{1}{r^{n-1}}\right)~\mbox{and}~|\nabla
G(0,y)|=O\left(\frac{1}{|ry|^{n-1}}\right).\ee Then, by
(\ref{eq-016csw}) and (\ref{eq-1c-2}),  there is a positive constant
$C_{1}$ such that

\beq\label{eqcxs-1} |\nabla
u(0)|&\leq&\frac{C_{1}}{r}\int_{\partial\mathbb{B}^{n}}|u(r\zeta)|d\sigma(\zeta)+\frac{C_{1}}{r}\sup_{\xi\in\mathbb{B}^{n}}\lambda(\xi)
\int_{\mathbb{B}^{n}}|u(ry)|^{\tau}|y|^{1-n}dy\\ \nonumber
&\leq&\frac{C_{1}Q_{R}(|u|)}{r}+\frac{nC_{1}V(\mathbb{B}^{n})\sup_{\xi\in\mathbb{B}^{n}}\lambda(\xi)Q_{R}(|u|^{\tau})}{r},
\eeq which, together with Lemma \ref{Lemx}, yield

\be\label{eq-1c-3}|\nabla
u(0)|^{\nu}\leq2^{\nu-1}C_{1}^{\nu}\left(\frac{Q_{R}(|u|^{\nu})}{r^{\nu}}+
\frac{n^{\nu}\sup_{\xi\in\mathbb{B}^{n}}\lambda^{\nu}(\xi)Q_{R}(|u|^{\nu\tau})V^{\nu}(\mathbb{B}^{n})}{r^{\nu}}\right),\ee
where
$Q_{R}(|u|)=\max\big\{|u(\xi)|:~\xi\in\overline{\mathbb{B}^{n}_{r}}\big\}$.

By (\ref{eq-1c-3}), Corollary \ref{cor-1} and Lemma \Ref{lem-eqc-3},
there is a positive constant $C_{2}$ such that

\begin{eqnarray*}|\nabla
u(0)|^{\nu}&\leq&2^{\nu-1}C_{1}^{\nu}C_{2}\bigg(\frac{1}{r^{\nu+n}}\int_{\mathbb{B}^{n}_{2r}}|u(y)|^{\nu}dy\\
&&+
\frac{n^{\nu}V^{\nu}(\mathbb{B}^{n})\sup_{\xi\in\mathbb{B}^{n}}\lambda^{\nu}(\xi)}{r^{n+\nu}}\int_{\mathbb{B}^{n}_{2r}}|u(y)|^{\tau\nu}dy\bigg).
\end{eqnarray*}
The proof of the proposition is complete. \qed

\begin{lem}\label{lem-03}
For  $\tau=1$, let $u\in \mathcal{C}^{2}(\mathbb{B}^{n})$ be a
solution to the equation {\rm (\ref{eq1c})}, where
 $\lambda$ is a nonnegative constant. For all
$a\in\mathbb{B}^{n}$, there is a positive constant $C$  such that
$$|\nabla u(a)|\leq
\frac{C}{r}\int_{\partial\mathbb{B}^{n}}|u(a+r\zeta)-u(a)|d\sigma(\zeta),$$
where $\overline{\mathbb{B}^{n}(a,r)}\subset\mathbb{B}^{n}$.
\end{lem}

\bpf For any fixed $a\in\mathbb{B}^{n}$, let
$f(x)=u(x+a)-u(a),~\mbox{$x\in\mathbb{B}^{n}_{r}$},$ where $r\in[0,
d(a)).$ By (\ref{eqcxs-1}) and  Corollary \ref{cor-1}, there is a
positive constant $C_{3}$ such that

\begin{eqnarray*} |\nabla
f(0)|&\leq&\frac{C_{3}}{r}\left(\int_{\partial\mathbb{B}^{n}}|f(r\zeta)|d\sigma(\zeta)+
\int_{\mathbb{B}^{n}}|f(ry)||y|^{1-n}dy\right)\\  &=&
\frac{C_{3}}{r}\left[\int_{\partial\mathbb{B}^{n}}|f(r\zeta)|d\sigma(\zeta)+n
\int_{0}^{1}\left(\int_{\partial\mathbb{B}^{n}}|f(r\rho\zeta)|d\sigma(\zeta)\right)d\rho\right]\\
&=&\frac{C_{3}}{r}\left[\int_{\partial\mathbb{B}^{n}}|f(r\zeta)|d\sigma(\zeta)+\frac{n}{r}
\int_{0}^{r}\Big(\int_{\partial\mathbb{B}^{n}}|f(t\zeta)|d\sigma(\zeta)\Big)dt\right]\\
&=&\frac{C_{3}}{r}\left[\int_{\partial\mathbb{B}^{n}}|f(r\zeta)|d\sigma(\zeta)+\frac{n}{r}
\int_{0}^{r}M_{1}(f,t)dt\right]\\
&\leq&\frac{C_{3}}{r}\left(\int_{\partial\mathbb{B}^{n}}|f(r\zeta)|d\sigma(\zeta)+nM_{1}(f,
r)\right)\\
&\leq&\frac{C_{3}(1+n)}{r}\int_{\partial\mathbb{B}^{n}}|f(r\zeta)|d\sigma(\zeta),
 \end{eqnarray*}
which yields $$|\nabla u(a)|\leq
\frac{C}{r}\int_{\partial\mathbb{B}^{n}}|u(a+r\zeta)-u(a)|d\sigma(\zeta),$$
where $C=C_{3}(1+n)$, completing the proof. \epf

\subsection*{Proof of  Theorem \ref{01-thm}}
First, we show the ``if" part. By Lemma \ref{lem-03}, there is a
positive constant $C_{4}$ such that

\be\label{eq-1c4}|\nabla u(x)|\leq
\frac{C_{4}}{\rho}\int_{\partial\mathbb{B}^{n}}|u(x+\rho\zeta)-u(x)|d\sigma(\zeta),\ee
where $\rho\in(0,d(x)].$ Let $r=d(x)$. Multiplying both sides of the
inequality (\ref{eq-1c4}) by $n\rho^{n-1}$ and integrating from $0$
to $r$, together with (\ref{eq-csw13}), we obtain

\begin{eqnarray*}
|\nabla
u(x)|&\leq&\frac{(n+1)C_{4}}{nr^{n+1}}\int_{0}^{r}\left(n\rho^{n-1}\int_{\partial\mathbb{B}^{n}}|u(x+\rho\zeta)-u(x)|d\sigma(\zeta)\right)d\rho\\
&=&\frac{(n+1)C_{4}}{nr|\mathbb{B}^{n}(x,r)|}\int_{\mathbb{B}^{n}(x,r)}|u(y)-u(x)|dy\\
&\leq&\frac{(n+1)C_{4}C}{n}\frac{1}{\omega(r^{\alpha})}\\
&=&\frac{(n+1)C_{4}C}{n}\frac{1}{\omega\big(d^{\alpha}(x)\big)},
\end{eqnarray*}
which implies that
$u\in\mathcal{L}_{+\infty,\omega}\mathcal{B}^{0}_{\alpha}(\mathbb{B}^{n}).$

Next we prove the ``only if" part. Since
$u\in\mathcal{L}_{+\infty,\omega}\mathcal{B}^{0}_{\alpha}(\mathbb{B}^{n}),$
we see that, for $x\in\mathbb{B}^{n}$, there is a positive constant
$C_{5}$ such that

\be\label{eq-1c-5}|\nabla
u(x)|\leq\frac{C_{5}}{\omega\big(d^{\alpha}(x)\big)}.\ee For $x,
y\in\mathbb{B}^{n}$ and $t\in[0,1]$, if $d(x)>t|x-y|$, then, by
(\ref{eq-1c-5}), we get

\begin{eqnarray*}
|u(x)-u(y)|&\leq&|x-y|\int_{0}^{1}|\nabla u(x+t(y-x))|dt\\
&\leq&C_{5}|x-y|\int_{0}^{1}\frac{dt}{\omega\big(d^{\alpha}(x+t(y-x))\big)}\\
&\leq&C_{5}|x-y|\int_{0}^{1}\frac{dt}{\omega\Big(\big(d(x)-t|x-y|\big)^{\alpha}\Big)}\\
&=&C_{5}\int_{0}^{|x-y|}\frac{dt}{\omega\left(\big(d(x)-t\big)^{\alpha}\right)},
\end{eqnarray*}
which yields that

\begin{eqnarray*}
I&\leq&
\frac{C_{5}}{|\mathbb{B}^{n}(x,r)|}\int_{\mathbb{B}^{n}(x,r)}\left[\int_{0}^{|x-y|}\frac{dt}{\omega\Big(\big(d(x)-t\big)^{\alpha}\Big)}
\right]dy\\
&=&\frac{C_{5}}{|\mathbb{B}^{n}_{r}|}\int_{\mathbb{B}^{n}_{r}}\left[\int_{0}^{|\xi|}\frac{dt}{\omega\Big(\big(d(x)-t\big)^{\alpha}\Big)}\right]d\xi\\
&=&\frac{C_{5}n}{r^{n}}\int_{0}^{r}\rho^{n-1}\left\{\int_{0}^{\rho}\frac{dt}{\omega\Big(\big(d(x)-t\big)^{\alpha}\Big)}
\right\}d\rho\\
&\leq&\frac{C_{5}n}{r^{n}}\int_{0}^{r}\left(\int_{t}^{r}\rho^{n-1}d\rho\right)\frac{1}{\omega\Big(\big(r-t\big)^{\alpha}\Big)}dt\\
&=&\frac{C_{5}}{r^{n}}\int_{0}^{r}\frac{(r-t)(r^{n-1}+r^{n-2}t+\cdots+t^{n-1})}{\omega\big((r-t)^{\alpha}\big)}dt\\
&\leq&\frac{C_{5}n}{r}\int_{0}^{r}\frac{(r-t)}{\omega\big((r-t)^{\alpha}\big)}dt\\
&=&\frac{C_{5}n}{r}\int_{0}^{r}\frac{(r-t)^{\alpha}}{\omega\big((r-t)^{\alpha}\big)}(r-t)^{1-\alpha}dt\\
&\leq&\frac{C_{5}nr^{\alpha-1}}{\omega(r^{\alpha})}\int_{0}^{r}(r-t)^{1-\alpha}dt\\
&=&\frac{C_{5}n}{(2-\alpha)}\frac{r}{\omega(r^{\alpha})},
\end{eqnarray*}
where
$$I=\frac{1}{|\mathbb{B}^{n}(x,r)|}\int_{\mathbb{B}^{n}(x,r)}|u(y)-u(x)|dy.$$
The proof of this theorem is complete.  \qed

For an $n\times n$ real  matrix $A$, we define the standard {\it operator
norm} by
$$\|A\|=\sup_{x\neq0}\frac{|Ax|}{|x|}=\max\big\{|A\theta|:\,\theta\in\partial\mathbb{B}^{n}\big\}.
$$

\subsection*{Proof of  Theorem \ref{02-thm}}
We first prove $(2)\Rightarrow(1)$. Let $x,y\in\mathbb{B}^{n}$ with
$r_{\mathbb{B}^{n}}(x,y)\leq1/2$. Then \be\label{eq-01csw}|x-y|\leq
d(x)/2.\ee By (\ref{eq-01csw}) and \cite[Lemma 3.7]{Vu}, we obtain
$$k_{\mathbb{B}^{n}}(x,y)\leq2j_{\mathbb{B}^{n}}(x,y)\leq2r_{\mathbb{B}^{n}}(x,y)\leq1,$$
which gives that
\be\label{eq-02csw}k_{u(\mathbb{B}^{n})}(u(x),u(y))\leq
Ck_{\mathbb{B}^{n}}(x,y)\leq C.\ee Applying (\ref{eq-02csw}), we get
$$j_{u(\mathbb{B}^{n})}(u(x),u(y))=\log\left(1+r_{u(\mathbb{B}^{n})}\big(u(x),u(y)\big)\right)\leq k_{u(\mathbb{B}^{n})}(u(x),u(y))\leq C,$$
which implies that $r_{u(\mathbb{B}^{n})}\big(u(x),u(y)\big)\leq
e^{C}-1$.

Now we prove  $(1)\Rightarrow(2)$. Since $u$ is weakly
uniformly bounded,  for every $x\in\mathbb{B}^{n}$ and
$y\in\overline{\mathbb{B}^{n}\left(x,d(x)/4\right)}$, we see  that
there is a positive constant $C$, \be\label{eq-03csw}|u(y)-u(x)|\leq
Cd_{u(\mathbb{B}^{n})}(u(x)).\ee By (\ref{eq-03csw}) and Lemma
\ref{lem-03}, we see that there is a positive $C_{6}$ such that

\beq\label{eq-06csw} \|u'(x)\|&\leq&\left(\sum_{k=1}^{n}|\nabla
u_{k}(x)|^{2}\right)^{\frac{1}{2}}\leq\sum_{k=1}^{n}|\nabla
u_{k}(x)|\\ \nonumber
&\leq&\frac{C_{6}}{r}\int_{\partial\mathbb{B}^{n}}\sum_{k=1}^{n}|u_{k}(x+r\zeta)-u_{k}(x)|d\sigma(\zeta)\\
\nonumber
&\leq&\frac{C_{6}\sqrt{n}}{r}\int_{\partial\mathbb{B}^{n}}|u(x+r\zeta)-u(x)|d\sigma(\zeta)\\
\nonumber
 &\leq&\frac{C_{6}C\sqrt{n}}{r}d_{u(\mathbb{B}^{n})}(u(x)), \eeq
where $r=d(x)$ and $$u'(x)=\left(\begin{array}{cccc}
\nabla u_{1}(x)   \\
\vdots \\
 \nabla u_{n}(x)
\end{array}\right).
$$
Hence $(1)\Rightarrow(2)$ follows from  (\ref{eq-06csw}) and
\cite[Lemma 2.6]{MM}.\qed

\begin{Thm}\label{Thm-B}
Let $g$ be a function of class $C^{2}(\mathbb{B}^{n})$. If $n\geq3$,
then for $r\in(0,1)$,
$$\int_{\partial \mathbb{B}^{n}}g(r\zeta)\,d\sigma(\zeta)=
g(0)+\int_{ \mathbb{B}^{n}(0,r)}\Delta g(x)G_{n}(x,r)\,dV_{N}(x),
$$
where $G_{n}(x,r)=(|x|^{2-n}-r^{2-n})/[n(n-2)]$ and $dV_{N}$ is the
normalized Lebesgue volume measure in $\mathbb{B}^{n}$. Moreover, if
$n=2$, then for $r\in(0,1)$,
\begin{eqnarray*}
\frac{1}{2\pi}\int_{0}^{2\pi}g(re^{i\theta})\,d\theta
&=&g(0)+ \frac{1}{2}\int_{\mathbb{D}_{r}}\Delta
g(z)\log\frac{r}{|z|}\,dA(z),
\end{eqnarray*}  where $dA$ denotes the normalized area measure in
$\mathbb{D}$  {\rm (cf. \cite{Pav,Z})}.
\end{Thm}

\begin{Lem}{\rm  (\cite[Lemma 3]{CPR})}\label{Lem-5}
Suppose that $\alpha>0$, $\beta\leq\alpha$ and $\omega$ is a
majorant. Then, for $r\in(0,1)$, $\phi(r)$  and
$\phi(r)/\omega(\phi(r))$ are decreasing in $(0,1)$, where $\phi$ is
the same as in Definition {\rm \ref{def-1}}.
\end{Lem}


\subsection*{Proof of  Theorem \ref{eq-y}}
 Without loss of generality, we assume that $u$ is a nonzero
function and $n\geq3$. By H\"older inequality, for $\rho\in[0,1),$
we have

$$M_{\frac{\nu(\nu-2)}{\nu-1}}^{\frac{\nu(\nu-2)}{\nu-1}}(u,\rho)\leq
\left(\int_{\partial\mathbb{B}^{n}}|u(\rho\zeta)|^{\nu}d\sigma(\zeta)\right)^{\frac{\nu-2}{\nu-1}}
\left(\int_{\partial\mathbb{B}^{n}}d\sigma(\zeta)\right)^{\frac{1}{\nu-1}},$$
which gives that

\beq\label{eq-y2}&&\int_{\partial\mathbb{B}^{n}}|u(\rho\zeta)|^{\nu-2}\big(|\nabla
u(\rho\zeta)|^{2}+u(\rho\zeta)\Delta
u(\rho\zeta)\big)d\sigma(\zeta)\\ \nonumber
&\leq&M_{\frac{\nu(\nu-2)}{\nu-1}}^{\nu-2}(u,\rho)\left[\int_{\partial\mathbb{B}^{n}}\Big(|\nabla
u(\rho\zeta)|^{2}+u(\rho\zeta)\Delta
u(\rho\zeta)\Big)^{\nu}d\sigma(\zeta)\right]^{\frac{1}{\nu}}\\
\nonumber
&\leq&M_{\nu}^{\nu-2}(u,\rho)\left[\int_{\partial\mathbb{B}^{n}}\Big(|\nabla
u(\rho\zeta)|^{2}+u(\rho\zeta)\Delta
u(\rho\zeta)\Big)^{\nu}d\sigma(\zeta)\right]^{\frac{1}{\nu}}. \eeq

By    (\ref{eq-y2}) and Theorem \Ref{Thm-B}, we obtain

 \begin{eqnarray*}
M_{\nu}^{\nu}(u,r)&=&|u(0)|^{\nu}+\int_{\mathbb{B}^{n}_{r}}\Delta(|u(x)|^{\nu})G_{n}(x,r)dV_{N}(x)\\
&=&|u(0)|^{\nu}+\int_{\mathbb{B}^{n}_{r}}\big[\nu(\nu-1)|u(x)|^{\nu-2}|\nabla
u(x)|^{2}\\
&&+\nu|u(x)|^{\nu-2}u(x)\Delta u(x)\big]G_{n}(x,r)dV_{N}(x)\\
&\leq&|u(0)|^{\nu}+\nu(\nu-1)\int_{0}^{r}\Big[n\rho^{n-1}G_{n}(\rho,r)\int_{\partial\mathbb{B}^{n}}|u(\rho\zeta)|^{\nu-2}\big(|\nabla
u(\rho\zeta)|^{2}\\
&&+u(\rho\zeta)\Delta u(\rho\zeta)\big)d\sigma(\zeta)\Big]d\rho\\
&\leq&|u(0)|^{\nu}+\nu(\nu-1)\int_{0}^{r}\bigg\{n\rho^{n-1}G_{n}(\rho,r)M_{\nu}^{\nu-2}(u,\rho)\\
&&\times\Big[\int_{\partial\mathbb{B}^{n}}\Big(|\nabla
u(\rho\zeta)|^{2}+u(\rho\zeta)\Delta
u(\rho\zeta)\Big)^{\nu}d\sigma(\zeta)\Big]^{\frac{1}{\nu}}\bigg\}d\rho\\
&\leq&|u(0)|^{\nu}+\nu(\nu-1)M_{\nu}^{\nu-2}(u,r)\int_{0}^{r}\bigg\{n\rho^{n-1}G_{n}(\rho,r)\\
&&\times\Big[\int_{\partial\mathbb{B}^{n}}\Big(|\nabla
u(\rho\zeta)|^{2}+u(\rho\zeta)\Delta
u(\rho\zeta)\Big)^{\nu}d\sigma(\zeta)\Big]^{\frac{1}{\nu}}\bigg\}d\rho,
\end{eqnarray*}
which, together with the subharmonicity of $u$ (Corollary
\ref{cor-1}) and Lemma \Ref{Lem-5}, yield that

\begin{eqnarray*}
M_{\nu}^{2}(u,r)&\leq&|u(0)|^{2}+\nu(\nu-1)\int_{0}^{r}\bigg\{n\rho^{n-1}G_{n}(\rho,r)\\
&&\times\Big[\int_{\partial\mathbb{B}^{n}}\Big(|\nabla
u(\rho\zeta)|^{2}+u(\rho\zeta)\Delta
u(\rho\zeta)\Big)^{\nu}d\sigma(\zeta)\Big]^{\frac{1}{\nu}}\bigg\}d\rho\\
&\leq&|u(0)|^{2}+\nu(\nu-1)C_{7}\int_{0}^{r}\frac{n\rho^{n-1}G_{n}(\rho,r)}{\omega\big(\phi(\rho)\big)}d\rho\\
&=&|u(0)|^{2}+\nu(\nu-1)C_{7}\int_{0}^{r}\frac{n\rho^{n-1}G_{n}(\rho,r)}{\phi(\rho)}\frac{\phi(\rho)}{\omega\big(\phi(\rho)\big)}d\rho\\
&\leq&|u(0)|^{2}+\frac{\nu(\nu-1)C_{7}}{\omega(1)}\int_{0}^{r}\frac{n\rho^{n-1}G_{n}(\rho,r)}{\phi(\rho)}d\rho\\
&=&|u(0)|^{2}+\frac{\nu(\nu-1)C_{7}r^{2}}{\omega(1)(n-2)}\int_{0}^{1}\frac{t(1-t^{n-2})}{\phi(tr)}dt,
\end{eqnarray*} where $C_{7}=\||\nabla u|^{2}+u\Delta
u\|_{\mathcal{L}_{\nu,\omega}\mathcal{B}^{\beta}_{\alpha}(\mathbb{B}^{n})}.$
The proof of the theorem is complete. \qed

\subsection*{Proof of  Theorem \ref{thm-1}}
 Without loss of generality, we assume that $u$ is a nonzero
function and $n\geq3$.  By H\"older inequality, for $\rho\in[0,1),$
we have
\be\label{eq-06}\int_{\partial\mathbb{B}^{n}}|u(\rho\zeta)|^{\nu-2}|\nabla
u(\rho\zeta)|^{2}d\sigma(\zeta)\leq
M_{\nu}^{\nu-2}(u,\rho)M_{\nu}^{2}(|\nabla u|, \rho), \ee

\be\label{eq-07}
\int_{\partial\mathbb{B}^{n}}|u(\rho\zeta)|^{\nu-1}|\nabla
u(\rho\zeta)|^{b_{1}}d\sigma(\zeta)\leq
M_{\nu}^{\nu-1}(u,\rho)M_{\nu b_{1}}^{b_{1}}(|\nabla u|, \rho),\ee

\be\label{eq-08} M_{\nu-1+b_{2}}^{\nu-1+b_{2}}( u,
\rho)\leq\left(\int_{\partial\mathbb{B}^{n}}|u(\rho\zeta)|^{\nu}d\sigma(\zeta)\right)^{\frac{\nu+b_{2}-1}{\nu}}
\left(\int_{\partial\mathbb{B}^{n}}d\sigma(\zeta)\right)^{\frac{1-b_{2}}{\nu}},
\ee

\be\label{eq-09} M_{\nu-1}^{\nu-1}( u,
\rho)\leq\left(\int_{\partial\mathbb{B}^{n}}|u(\rho\zeta)|^{\nu}d\sigma(\zeta)\right)^{\frac{\nu-1}{\nu}}
\left(\int_{\partial\mathbb{B}^{n}}d\sigma(\zeta)\right)^{\frac{1}{\nu}},
\ee and

\be\label{eq-10} M_{\nu b_{1}}^{\nu b_{1}}( |\nabla u|,
\rho)\leq\left(\int_{\partial\mathbb{B}^{n}}|\nabla
u(\rho\zeta)|^{\nu}d\sigma(\zeta)\right)^{\frac{\nu b_{1}}{\nu}}
\left(\int_{\partial\mathbb{B}^{n}}d\sigma(\zeta)\right)^{\frac{\nu-\nu
b_{1}}{\nu}}. \ee

Applying   (\ref{eq-06}),  (\ref{eq-07}), (\ref{eq-08}),
(\ref{eq-09}), (\ref{eq-10}), \cite[Lemma 3]{CPR} and Theorem
\Ref{Thm-B}, for $r\in[0,1)$, we get
\begin{eqnarray*}
M_{\nu}^{\nu}(u,r)&=&|u(0)|^{\nu}+\int_{\mathbb{B}^{n}_{r}}\Delta(|u(x)|^{\nu})G_{n}(x,r)dV_{N}(x)\\
&=&\nu(\nu-1)\int_{0}^{r}\left[n\rho^{n-1}G_{n}(\rho,r)\int_{\partial\mathbb{B}^{n}}|u(\rho\zeta)|^{\nu-2}|\nabla
u(\rho\zeta)|^{2}d\sigma(\zeta)\right]d\rho\\
&&+\nu\int_{0}^{r}\left[n\rho^{n-1}G_{n}(\rho,r)\int_{\partial\mathbb{B}^{n}}|u(\rho\zeta)|^{\nu-2}u(\rho\zeta)\Delta u(\rho\zeta)d\sigma(\zeta)\right]d\rho\\
&&+|u(0)|^{\nu}\\
&\leq&|u(0)|^{\nu}+\nu(\nu-1)\int_{0}^{r}n\rho^{n-1}G_{n}(\rho,r)M_{\nu}^{\nu-2}(u,\rho)M_{\nu}^{2}(|\nabla u|, \rho)d\rho\\
&&+\nu\sup_{x\in\mathbb{B}^{n}}a_{1}(x)\int_{0}^{r}n\rho^{n-1}G_{n}(\rho,r)M_{\nu}^{\nu-1}(u,\rho)M_{\nu
b_{1}}^{b_{1}}(|\nabla u|, \rho)d\rho\\
&&+\nu\sup_{x\in\mathbb{B}^{n}}a_{2}(x)\int_{0}^{r}n\rho^{n-1}G_{n}(\rho,r)M_{\nu-1+b_{2}}^{\nu-1+b_{2}}(
u, \rho)d\rho\\
&&+\nu\sup_{x\in\mathbb{B}^{n}}a_{3}(x)\int_{0}^{r}n\rho^{n-1}G_{n}(\rho,r)M_{\nu-1}^{\nu-1}(
u, \rho)d\rho\\
&\leq&|u(0)|^{\nu}+\nu(\nu-1)\int_{0}^{r}n\rho^{n-1}G_{n}(\rho,r)M_{\nu}^{\nu-2}(u,\rho)M_{\nu}^{2}(|\nabla u|, \rho)d\rho\\
&&+\nu\sup_{x\in\mathbb{B}^{n}}a_{1}(x)\int_{0}^{r}n\rho^{n-1}G_{n}(\rho,r)M_{\nu}^{\nu-1}(u,\rho)M_{\nu
b_{1}}^{b_{1}}(|\nabla u|, \rho)d\rho\\
&&+\nu\sup_{x\in\mathbb{B}^{n}}a_{2}(x)\int_{0}^{r}n\rho^{n-1}G_{n}(\rho,r)M_{\nu}^{\nu-1+b_{2}}(
u, \rho)d\rho\\
&&+\nu\sup_{x\in\mathbb{B}^{n}}a_{3}(x)\int_{0}^{r}n\rho^{n-1}G_{n}(\rho,r)M_{\nu}^{\nu-1}(
u, \rho)d\rho,
\end{eqnarray*}
which gives that

\begin{eqnarray*}
M_{\nu}^{2}(u,r)&\leq&|u(0)|^{2}+\nu(\nu-1)\int_{0}^{r}n\rho^{n-1}G_{n}(\rho,r)M_{\nu}^{2}(|\nabla u|, \rho)d\rho\\
&&+\nu
M_{\nu}(u,r)\sup_{x\in\mathbb{B}^{n}}a_{1}(x)\int_{0}^{r}n\rho^{n-1}G_{n}(\rho,r)M_{\nu
b_{1}}^{b_{1}}(|\nabla u|, \rho)d\rho\\
&&+\nu M_{\nu}^{1+b_{2}}(
u, r)\sup_{x\in\mathbb{B}^{n}}a_{2}(x)\int_{0}^{r}n\rho^{n-1}G_{n}(\rho,r)d\rho\\
&&+\nu M_{\nu}( u,
r)\sup_{x\in\mathbb{B}^{n}}a_{3}(x)\int_{0}^{r}n\rho^{n-1}G_{n}(\rho,r)d\rho\\
\end{eqnarray*}
\begin{eqnarray*}
&=&|u(0)|^{2}+\frac{\nu(\nu-1)r^{2}}{n-2}\int_{0}^{1}t(1-t^{n-2})M_{\nu}^{2}(|\nabla u|, rt)dt\\
&&+\frac{\nu
r^{2}\sup_{x\in\mathbb{B}^{n}}a_{1}(x)}{n-2}M_{\nu}(u,r)\int_{0}^{1}t(1-t^{n-2})M_{\nu
}^{b_{1}}(|\nabla u|, rt)dt\\
&&+\frac{\nu
r^{2}\sup_{x\in\mathbb{B}^{n}}a_{2}(x)}{2n}M_{\nu}^{1+b_{2}}(
u, r)\\
&&+\frac{\nu r^{2}\sup_{x\in\mathbb{B}^{n}}a_{3}(x)}{2n}M_{\nu}( u,
r)\\
&\leq&|u(0)|^{2}+\frac{\nu(\nu-1)}{n-2}\|u\|^{2}_{\mathcal{L}_{\nu,\omega}
\mathcal{B}^{\beta}_{\alpha}(\mathbb{B}^{n})}r^{2}\int_{0}^{1}\frac{\phi^{2}(rt)}{\omega^{2}\big(\phi(rt)\big)}\frac{t(1-t^{n-2})}{\phi^{2}(rt)}dt\\
&&+\frac{\nu
\sup_{x\in\mathbb{B}^{n}}a_{1}(x)}{n-2}\|u\|^{b_{1}}_{\mathcal{L}_{\nu,\omega}
\mathcal{B}^{\beta}_{\alpha}(\mathbb{B}^{n})}r^{2}M_{\nu}(u,r)\int_{0}^{1}\frac{\phi^{b_{1}}(rt)}{\omega^{b_{1}}\big(\phi(rt)\big)}\frac{t(1-t^{n-2})}{\phi^{b_{1}}(rt)}dt\\
&&+\frac{\nu
\sup_{x\in\mathbb{B}^{n}}a_{2}(x)}{2n}r^{2}M_{\nu}^{1+b_{2}}(
u, r)\\
&&+\frac{\nu \sup_{x\in\mathbb{B}^{n}}a_{3}(x)}{2n}r^{2}M_{\nu}( u,
r)\\
&\leq&|u(0)|^{2}+\frac{\nu(\nu-1)}{(n-2)\omega^{2}(1)}\|u\|^{2}_{\mathcal{L}_{\nu,\omega}
\mathcal{B}^{\beta}_{\alpha}(\mathbb{B}^{n})}r^{2}\int_{0}^{1}\frac{t(1-t^{n-2})}{\phi^{2}(rt)}dt\\
&&+\frac{\nu
\sup_{x\in\mathbb{B}^{n}}a_{1}(x)}{(n-2)\omega^{b_{1}}(1)}\|u\|^{b_{1}}_{\mathcal{L}_{\nu,\omega}
\mathcal{B}^{\beta}_{\alpha}(\mathbb{B}^{n})}r^{2}M_{\nu}(u,r)\int_{0}^{1}\frac{t(1-t^{n-2})}{\phi^{b_{1}}(rt)}dt\\
&&+\frac{\nu
\sup_{x\in\mathbb{B}^{n}}a_{2}(x)}{2n}r^{2}M_{\nu}^{1+b_{2}}(
u, r)\\
&&+\frac{\nu \sup_{x\in\mathbb{B}^{n}}a_{3}(x)}{2n}r^{2}M_{\nu}( u,
r).
\end{eqnarray*}
The proof of the theorem is complete.
 \qed

\begin{lem}\label{lem-cw4}
 Let  $u\in \mathcal{C}^{2}(\mathbb{B}^{n})$ be a solution to
the equation {\rm (\ref{eq1c})} with $\tau=1$. Then, for $\nu\geq1$,
$\left(\sum_{1\leq k,j\leq n}u^{2}_{x_{k}x_{j}}\right)^{\nu}$ is
subharmonic in $\mathbb{B}^{n}$.
\end{lem}

\bpf Without loss of generality, we assume that $U=\left(\sum_{1\leq
k,j\leq n}u^{2}_{x_{k}x_{j}}\right)^{\nu}$
has no zeros. 
By computations, we get

\begin{eqnarray*}
\Delta U&=&4\nu(\nu-1)\left(\sum_{1\leq k,j\leq
n}u^{2}_{x_{k}x_{j}}\right)^{\nu-2}
\left(\sum_{1\leq m,k,j\leq n}u_{x_{k}x_{j}}u_{x_{k}x_{j}x_{m}}\right)^{2}\\
&&+2\nu\left(\sum_{1\leq k,j\leq n}u^{2}_{x_{k}x_{j}}\right)^{\nu-1}
\sum_{1\leq k,j\leq n}\left[(\Delta
u)_{x_{k}x_{j}}u_{x_{k}x_{j}}+\sum_{m=1}^{n}u^{2}_{x_{m}x_{k}x_{j}}\right]\\
\end{eqnarray*}
\begin{eqnarray*}
&=&4\nu(\nu-1)\left(\sum_{1\leq k,j\leq
n}u^{2}_{x_{k}x_{j}}\right)^{\nu-2}
\left(\sum_{1\leq m,k,j\leq n}u_{x_{k}x_{j}}u_{x_{k}x_{j}x_{m}}\right)^{2}\\
&&+2\nu\left(\sum_{1\leq k,j\leq n}u^{2}_{x_{k}x_{j}}\right)^{\nu-1}
\sum_{1\leq k,j\leq n}\left[\lambda
u^{2}_{x_{k}x_{j}}+\sum_{m=1}^{n}u^{2}_{x_{m}x_{k}x_{j}}\right]\geq0,
\end{eqnarray*}
which implies that $U$ is  subharmonic in $\mathbb{B}^{n}$.  \epf

\subsection*{Proof of  Theorem \ref{thm-3}} By Lemma \Ref{lem-eqc-3} and Lemma \ref{lem-cw4}, for
$\mu\in[1,n/2]$ and $x\in\mathbb{B}^{n}$, there is a positive
constant $C_{8}$ such that

$$\left(\sum_{1\leq k,j\leq
n}u^{2}_{x_{k}x_{j}}(x)\right)^{\mu}\leq\frac{C_{8}}{(d(x))^{n+\alpha}}\int_{\mathbb{B}^{n}\left(x,\frac{d(x)}{2}\right)}(1-|y|)^{\alpha}\left(\sum_{1\leq
k,j\leq n}u^{2}_{x_{k}x_{j}}(y)\right)^{\mu}dy,$$ which gives that

\be\label{eq-csw5} \sum_{1\leq k,j\leq
n}u^{2}_{x_{k}x_{j}}(x)\leq\frac{\big(C_{8}\mathcal{D}_{\nabla
u}(\alpha,0,\mu)\big)^{\frac{1}{\mu}}}{(d(x))^{\frac{n+\alpha}{\mu}}}.
\ee

Let $$H_{u}=\left(\begin{array}{cccc} \ds \frac{\partial^{2}
u}{\partial x^{2}_{1}}\; \frac{\partial^{2} u}{\partial
x_{1}\partial x_{2}}\; \cdots\;
  \frac{\partial^{2} u}{\partial x_{1}\partial x_{n}}\\[4mm]
 \ds  \frac{\partial^{2} u}{\partial x_{2}\partial x_{1}}\; \frac{\partial^{2}
u}{\partial x^{2}_{2}}\; \cdots\;
  \frac{\partial^{2} u}{\partial x_{2}\partial x_{n}}\\[2mm]
\vdots \\[2mm]
\ds \frac{\partial^{2} u}{\partial x_{n}\partial x_{1}}\;
\frac{\partial^{2} u}{\partial x_{n}\partial x_{2}}\; \cdots\;
  \frac{\partial^{2} u}{\partial x_{n}^{2}}
\end{array}\right)
$$ be the  Hessian matrix of $u$. Then

\be\label{eq-csw0} \|H_{u}\|\leq\sqrt{\sum_{1\leq k,j\leq
n}u^{2}_{x_{k}x_{j}}}.\ee

By (\ref{eq-csw5}) and (\ref{eq-csw0}), we get

\beq\label{eq-csw6} |\nabla u(x)|&\leq& |\nabla u(0)|+
\int_{[0,x]}\|H_{u}(y)\||dy|\\ \nonumber&\leq&|\nabla
u(0)|+\int_{[0,x]}\frac{\big(C_{8}\mathcal{D}_{\nabla
u}(\alpha,0,\mu)\big)^{\frac{1}{2\mu}}}{(d(y))^{\frac{n+\alpha}{2\mu}}}|dy|\\
\nonumber &\leq&|\nabla
u(0)|+\frac{C_{9}}{(d(x))^{\frac{n+\alpha}{2\mu}-1}}, \eeq
where
\[
C_{9}=\frac{2\mu\big(C_{8}\mathcal{D}_{\nabla
u}(\alpha,0,\mu)\big)^{\frac{1}{2\mu}}}{n+\alpha-2\mu},
\]  and $[0,x]$ is the line segment from $0$ to $x$.



Applying (\ref{eq-csw6}) and Lemma \ref{Lemx}, for $\nu\geq2$, we
have

\beq\label{eq-csw7} |\nabla u(x)|^{\nu-2}&\leq&\left[|\nabla
u(0)|+\frac{C_{9}}{(d(x))^{\beta}}\right]^{\nu-2}\\
\nonumber &\leq&2^{\nu-2}\left[|\nabla
u(0)|^{\nu-2}+\frac{C_{9}^{\nu-2}}{(d(x))^{\beta(\nu-2)}}\right]\eeq

and

\beq\label{eq-csw8} |\nabla u(x)|^{\nu}&\leq&\left[|\nabla
u(0)|+\frac{C_{9}}{(d(x))^{\beta}}\right]^{\nu}\\
\nonumber &\leq&2^{\nu}\left[|\nabla
u(0)|^{\nu}+\frac{C_{9}^{\nu}}{(d(x))^{\beta\nu}}\right],\eeq where
$\beta=\frac{n+\alpha}{2\mu}-1.$

We divide the remaining part of the proof into two cases, namely
$\nu\in[2,4)$ and $\nu\in[4,+\infty).$

$\mathbf{Case~ I:}$ Let $ \nu\in[4,+\infty).$  By direct
computations, we see that

\beq\label{eq-csw9} \Delta \left(|\nabla u|^{\nu}\right)
&=&\nu(\nu-2)|\nabla
u|^{\nu-4}\sum_{j=1}^{n}\left(\sum_{k=1}^{n}u_{x_{k}x_{j}}u_{x_{k}}\right)^{2}\\
\nonumber &&+\nu|\nabla u|^{\nu-2}\sum_{k=1}^{n}u_{x_{k}}(\Delta
u)_{x_{k}}+\nu|\nabla
u|^{\nu-2}\sum_{j=1}^{n}\sum_{k=1}^{n}u^{2}_{x_{k}x_{j}}\\ \nonumber
&\leq&\nu(\nu-2)|\nabla u|^{\nu-2}\sum_{1\leq k,j\leq
n}u^{2}_{x_{k}x_{j}}+\lambda \nu|\nabla u|^{\nu}\\
\nonumber &&+\nu|\nabla
u|^{\nu-2}\sum_{1\leq k,j\leq n}u^{2}_{x_{k}x_{j}}\\
\nonumber&=&\nu(\nu-1)|\nabla u|^{\nu-2}\sum_{1\leq k,j\leq
n}u^{2}_{x_{k}x_{j}}+\lambda \nu|\nabla u|^{\nu}. \eeq

It follows from (\ref{eq-csw7}), (\ref{eq-csw8}) and (\ref{eq-csw9})
that

\beq\label{eq-csw10}  \big(d(x)\big)^{\beta\nu}\Delta \left(|\nabla
u|^{\nu}\right)&\leq&\nu(\nu-1)\big(d(x)\big)^{\beta\nu}|\nabla
u|^{\nu-2}\sum_{1\leq k,j\leq n}u^{2}_{x_{k}x_{j}}\\ \nonumber
&&+\lambda \nu\big(d(x)\big)^{\beta\nu}|\nabla u|^{\nu}\\
\nonumber&=&\nu(\nu-1)\big(d(x)\big)^{\beta\nu-\frac{\alpha}{\mu}}|\nabla
u|^{\nu-2}\big(d(x)\big)^{\frac{\alpha}{\mu}}\sum_{1\leq k,j\leq n}u^{2}_{x_{k}x_{j}}\\
\nonumber &&+\lambda \nu\big(d(x)\big)^{\beta\nu}|\nabla u|^{\nu}\\
\nonumber &\leq&C_{10}\big(d(x)\big)^{\frac{\alpha}{\mu}}\sum_{1\leq
k,j\leq n}u^{2}_{x_{k}x_{j}}+C_{11}, \eeq where
$C_{10}=2^{\nu-2}\nu(\nu-1)\big(|\nabla
u(0)|^{\nu-2}+C_{9}^{\nu-2}\big)$ and $C_{11}=2^{\nu}\lambda
\nu\big(|\nabla u(0)|^{\nu}+C_{9}^{\nu}\big).$

By H\"older's inequality, we obtain

\beq\label{eq-csw11}
\int_{\mathbb{B}^{n}}\big(d(x)\big)^{\frac{\alpha}{\mu}}\left(\sum_{1\leq
k,j\leq n}u^{2}_{x_{k}x_{j}}(x)\right)dx&\leq&
\big(\mathcal{D}_{\nabla
u}(\alpha,0,\mu)\big)^{\frac{1}{\mu}}\left(\int_{\mathbb{B}^{n}}dx\right)^{1-\frac{1}{\mu}}\\
\nonumber &=&
\big(V(\mathbb{B}^{n})\big)^{1-\frac{1}{\mu}}\big(\mathcal{D}_{\nabla
u}(\alpha,0,\mu)\big)^{\frac{1}{\mu}}. \eeq

By (\ref{eq-csw10}) and (\ref{eq-csw11}), we conclude that

\beq\label{eq-csw12} \nonumber
\int_{\mathbb{B}^{n}}\big(d(x)\big)^{\beta\nu}\Delta \left(|\nabla
u(x)|^{\nu}\right)dx&\leq& \int_{\mathbb{B}^{n}}\left[
C_{10}\big(d(x)\big)^{\frac{\alpha}{\mu}}\sum_{1\leq k,j\leq
n}u^{2}_{x_{k}x_{j}}(x)+C_{11}\right]dx\\
&\leq&C_{10}\big(V(\mathbb{B}^{n})\big)^{2-\frac{1}{\mu}}\big(\mathcal{D}_{\nabla
u}(\alpha,0,\mu)\big)^{\frac{1}{\mu}}+C_{11}V(\mathbb{B}^{n})\\
\nonumber &<&+\infty. \eeq

$\mathbf{Case~ II:}$ Let $ \nu\in[2,4).$ In this case, for
$m\in\{1,2,\ldots\}$, we let $f_{m}^{\nu}=\left(|\nabla
u|^{2}+\frac{1}{m}\right)^{\frac{\nu}{2}}$. It is not difficult to
know that $\Delta(f_{m}^{\nu})$ is integrable in
$\mathbb{B}^{n}_{r}$. Then, by (\ref{eq-csw10}), (\ref{eq-csw12})
and Lebesgue's dominated convergence theorem, we have

\begin{eqnarray*}
\lim_{m\rightarrow+\infty}\int_{\mathbb{B}^{n}}\big(d(x)\big)^{\beta\nu}\Delta
\left(
f_{m}^{\nu}(x)\right)dx&=&\int_{\mathbb{B}^{n}}\big(d(x)\big)^{\beta\nu}\lim_{m\rightarrow+\infty}\Delta
\left( f_{m}^{\nu}(x)\right)dx\\
&\leq&\int_{\mathbb{B}^{n}}\left[
C_{10}\big(d(x)\big)^{\frac{\alpha}{\mu}}\sum_{1\leq k,j\leq
n}u^{2}_{x_{k}x_{j}}(x)+C_{11}\right]dx\\
&\leq&C_{10}\big(V(\mathbb{B}^{n})\big)^{2-\frac{1}{\mu}}\big(\mathcal{D}_{\nabla
u}(\alpha,0,\mu)\big)^{\frac{1}{\mu}}+C_{11}V(\mathbb{B}^{n})\\
 &<&+\infty.
\end{eqnarray*}
The proof of the theorem is complete. \qed

\begin{lem}\label{lem-cw5}
Let $u\in \mathcal{C}^{3}(\mathbb{B}^{n})$ with
$\sum_{k=1}^{n}u_{x_{k}}(\Delta u)_{x_{k}}\geq0$ in
$\mathbb{B}^{n}$. Then, for $\nu\geq1$, $|\nabla u|^{\nu}$ is
subharmonic in $\mathbb{B}^{n}$.

\end{lem}

\bpf Let $\mathcal{Z}_{\nabla u}=\{x\in\mathbb{B}^{n}:~|\nabla
u(x)|=0\}$. Then  $\mathbb{B}^{n}\backslash\mathcal{Z}_{\nabla u}$
is an open set. For $j\in\{1,\ldots,n\}$ and
$x\in\mathbb{B}^{n}\backslash\mathcal{Z}_{\nabla u}$, we have

\begin{eqnarray*}
(|\nabla u(x)|^{\nu})_{x_{j}x_{j}}&=&\nu(\nu-2)|\nabla
u(x)|^{\nu-4}\left(\sum_{k=1}^{n}u_{x_{k}x_{j}}(x)u_{x_{k}}(x)\right)^{2}\\
&&+\nu|\nabla
u(x)|^{\nu-2}\sum_{k=1}^{n}\left(u_{x_{k}x_{j}x_{j}}(x)u_{x_{k}}(x)+u_{x_{k}x_{j}}^{2}(x)\right),
\end{eqnarray*}
which gives that

\begin{eqnarray*}
\Delta \left(|\nabla u(x)|^{\nu}\right) &=&\nu(\nu-2)|\nabla
u(x)|^{\nu-4}\sum_{j=1}^{n}\left(\sum_{k=1}^{n}u_{x_{k}x_{j}}(x)u_{x_{k}}(x)\right)^{2}\\
&&+\nu|\nabla u(x)|^{\nu-2}\sum_{k=1}^{n}u_{x_{k}}(x)(\Delta
u(x))_{x_{k}}+\nu|\nabla
u(x)|^{\nu-2}\sum_{j=1}^{n}\sum_{k=1}^{n}u^{2}_{x_{k}x_{j}}(x)\\
&\geq&0.
\end{eqnarray*}
Therefore,  for $\nu\geq1$, $|\nabla u|^{\nu}$ is subharmonic in
$\mathbb{B}^{n}$. \epf

The following result easily follows from Lemma \ref{lem-cw5}.

\begin{cor}\label{cor-cxs2}
Let $u\in \mathcal{C}^{3}(\mathbb{B}^{n})$ be a solution to the
equation {\rm (\ref{eq1c})}, where $\lambda$ is a nonnegative
constant. Then, for $\nu\geq1$, $|\nabla u|^{\nu}$ is subharmonic in
$\mathbb{B}^{n}$.
\end{cor}

\subsection*{Proof of  Theorem \ref{thm-4}} $|\nabla
u|\in\mathcal{H}^{\nu}_{g}(\mathbb{B}^{n})$  follows from
\cite[Theorem 1]{CRV} and Theorem \ref{thm-3}.

Next we prove $|\nabla u|^{\nu}$ has a harmonic majorant. For
$x\in\mathbb{B}^{n}$, let
$$G_{r}(x)=\int_{\partial\mathbb{B}^{n}}\frac{1-|x|^{2}}{|x-\zeta|^{n}}|\nabla
u(r\zeta)|^{\nu}d\sigma(\zeta),$$ where $r\in[0,1).$ By Corollary
\ref{cor-cxs2}, we see that $|\nabla u|^{\nu}$ is subharmonic in
$\mathbb{B}^{n}$, which, together with $|\nabla
u|\in\mathcal{H}^{\nu}_{g}(\mathbb{B}^{n})$, imply that

$$|\nabla u(rx)|^{\nu}\leq\int_{\partial\mathbb{B}^{n}}\frac{1-|x|^{2}}
{|x-\zeta|^{n}}|\nabla
u(r\zeta)|^{\nu}d\sigma(\zeta)=G_{r}(x)<+\infty$$  and
$G_{r}(0)=M_{\nu}^{\nu}(|\nabla u|,r)<+\infty$. For
$x\in\mathbb{B}^{n}$, applying the Harnack Theorem  to the sequence
$\{G_{1-1/m}(x)\}_{m=1}^{\infty}$, we see that
$$G(x)=\lim_{m\rightarrow+\infty}G_{1-1/m}(x)$$ is also a harmonic
function in $\mathbb{B}^{n}$. Hence $|\nabla u|^{\nu}$ has a
harmonic majorant in $\mathbb{B}^{n}$. The proof of the theorem is
complete. \qed



\normalsize

\end{document}